\documentclass[a4paper,11pt]{amsart}
\usepackage{graphicx}
\usepackage{mathptmx}
\usepackage{amsmath}
\usepackage{amssymb}
\usepackage{enumitem}
\usepackage{xcolor}
\usepackage{xparse}
\NewDocumentCommand{\eulerian}{omm}
 {%
  \genfrac<>{0pt}{}{#2}{#3}%
  \IfValueT{#1}{_{\!#1}}%
 }

 \newmuskip\pFqmuskip

\newcommand*\pFq[6][8]{%
  \begingroup % only local assignments
  \pFqmuskip=#1mu\relax
  \mathchardef\normalcomma=\mathcode`,
  \mathcode`\,=\string"8000
  \begingroup\lccode`\~=`\,
  \lowercase{\endgroup\let~}\pFqcomma
  {}_{#2}F_{#3}{\left(\genfrac..{0pt}{}{#4}{#5}\bigg|#6\right)}%
  \endgroup
}
\newcommand{\pFqcomma}{{\normalcomma}\mskip\pFqmuskip}

\newtheorem{theorem}{Theorem}

\begin{document}

\title[degenerate $r$-Whitney numbers and degenerate $r$-Dowling polynomials]
{Degenerate $r$-Whitney numbers and degenerate $r$-Dowling polynomials via boson operators}

\author{Taekyun  Kim}
\address{Department of Mathematics, Kwangwoon University, Seoul 139-701, Republic of Korea}
\email{tkkim@kw.ac.kr}

\author{DAE SAN KIM}
\address{Department of Mathematics, Sogang University, Seoul 121-742, Republic of Korea}
\email{dskim@sogang.ac.kr}

\subjclass[2010]{11B73; 11B83}
\keywords{degenerate $r$-Whitney numbers of the second kind; degenerate $r$-Whitney numbers of the first kind; degenerate $r$-Dowling polynomials; boson operator; normal ordering}

\maketitle

\begin{abstract}
Dowling showed that the Whitney numbers of the first kind and of the second kind satisfy Stirling number-like relations. Recently, Kim-Kim introduced the degenerate $r$-Whitney numbers of the first kind and of the second kind, as degenerate versions and further generalizations of the Whitney numbers of both kinds. The normal ordering of an integral power of the number operator in terms of boson operators is expressed with the help of the Stirling numbers of the second kind. In this paper, it is noted that the normal ordering of a certain quantity involving the number operator is expressed in terms of the degenerate $r$-Whitney numbers of the second kind. We derive some properties, recurrence relations, orthogonality relations and several identities on those numbers from such normal ordering.
In addition, we consider the degenerate $r$-Dowling polynomials as a natural extension of the degenerate $r$-Whitney numbers of the second kind and investigate their properties. 
\end{abstract}

\section{Introduction}
Dowling [3,4] constructed an important finite geometric lattice $Q_{n}(G)$ out of a finite set of $n$ elements and a finite group $G$ of order $m$, called Dowling lattice of rank $n$ over a finite group of order $m$. For a finite geometric lattice $L$ of rank $n$, Dowling [3,4] defined the Whitney numbers of the first kind $V_{L}(n,k)$ and the Whitney numbers of the second $W_{L}(n,k)$. If $L$ is the Dowling lattice $Q_{n}(G)$ of rank $n$ over finite group $G$ of order $m$. then the Whitney numbers of the first kind $V_{Q_{n}(G)}(n,k)$ and the Whitney numbers of the second kind $W_{Q_{n}(G)}(n,k)$ are respectively denoted by $V_{m}(n,k)$ and $W_{m}(n,k)$. These notations are justified, since the Whitney numbers of both kinds depend only on the order $m$ of $G$. It is known that, for any fixed group $G$ of order $m$, the Whitney numbers $W_{m}(n,k)$ and $V_{m}(n,k)$ satisfy the Stirling number-like relations in \eqref{8} and \eqref{9}. Recently, the study of degenerate versions of some special numbers and polynomials has regained interests of quite a few mathematicians. Kim-Kim introduced in [10] the degenerate $r$-Whitney numbers of the first kind $V_{m,\lambda}^{(r)}(n,k)$ and of the second $W_{m,\lambda}^{(r)}(n,k)$, as degenerate versions and extensions of the Whitney numbers of both kinds.They satisfy the relations in \eqref{15} and \eqref{16}. \par
The normal ordering of an integral power of the number operator $a^{\dagger}a$ in terms of boson operators $a$ and $a^{\dagger}$ can be written in the form 
\begin{equation*}
(a^{+}a)^{k}=\sum_{l=0}^{k}S_{2}(k,l)(a^{+})^{l}a^{l},
\end{equation*}
and, by inversion, we also have
\begin{equation*}
(a^{+})^{k}a^{k}=\sum_{l=0}^{k}S_{1}(k,l)(a^{+}a)^{l}.
\end{equation*}
The normal ordering of $(ma^{+}a+r)_{k,\lambda}$ in terms of boson operators is given in the form
\begin{equation}
(ma^{+}a+r)_{k,\lambda}=\sum_{l=0}^{k}W_{m,\lambda}^{(r)}(k,l)m^{l} (a^{+})^{l}a^{l}, \label{-2}
\end{equation}
and, by inversion, we also have
\begin{equation}
m^{k}(a^{+})^{k}a^{k}=\sum_{l=0}^{k}V_{m,\lambda}^{(r)}(k,l)(ma^{+}a+r)_{l,\lambda}. \label{-1}
\end{equation}
The aim of this paper is to derive some properties, recurrence relations, orthogonality relations and several identities on the degenerate $r$-Whitney numbers of both kinds from the normal ordering in \eqref{-2} and its inversion in \eqref{-1}.
In addition, we consider the degenerate $r$-Dowling polynomials as a natural extension of the degenerate $r$-Whitney numbers of the second kind and investigate their properties. The novelty of this paper is that the degenerate $r$-Whitney numbers are investigated by boson operators, the number operator, normal ordering and coherent states in quantum mechanics. 
Here we mention two recent works written in the same spirit as the present paper. In [12], it is noted that the normal ordering of a `degenerate integral power' of the number operator in terms of boson operators is represented by means of the degenerate Stirling numbers of the second kind. As an application of this normal ordering, obtained are two equations defining the degenerate Stirling numbers of the second kind and a Dobinski-like formula for the degenerate Bell polynomials. In [9], deduced are some identities and recurrence relations for the degenerate Stirling numbers of the first kind and of the second kind from the normal oderings of the degenerate integral powers of the number operator and their inversions and certain relations of boson operators as well. \par
The outline of this paper is as follows. In Section 1, we recall some definitions, including the degenerate exponential functions, the degenerate Stirling numbers of the first kind $S_{1,\lambda}(n,k)$ and of the second kind $S_{2,\lambda}(n,k)$, the unsigned degenerate $r$-Stirling numbers of the first kind ${n+r \brack k+r}_{r,\lambda}$ and the degenerate $r$-Stirling numbers of the second kind ${n+r \brace k+r}_{r,\lambda}$. And then we recall the Whitney numbers of the first kind $V_{m}(n,k)$ and of the second kind $W_{m}(n,k)$, the $r$-Whitney numbers of the first kind $V_{m}^{(r)}(n,k)$ and of the second kind $W_{m}^{(r)}(n,k)$, the degenerate Whitney numbers of the first kind $V_{m,\lambda}(n,k)$ and of the second kind $W_{m,\lambda}(n,k)$, the degenerate $r$-Whitney numbers of the first kind $V_{m,\lambda}^{(r)}(n,k)$ and of the second kind $W_{m,\lambda}^{(r)}(n,k)$ and the degenerate Dowling polynomials $D_{m,\lambda}(n,x)$.                          
And then we recall the commutation relation for boson operators $a$ and  $a^{+}$ (see \eqref{18}), the normal ordering of an integral power of the number operator $a^{+}a$ (see \eqref{19}) and the number states (see \eqref{20}).
Section 2 is the main results of this paper. In Section 2, we show that the normal ordering of $(ma^{+}a+r)_{n,\lambda}$ is expressed in terms of the degenerate $r$-Whitney numbers of the second kind so that its inversion is represented by the degenerate $r$-Whitney numbers of the first kind (see \eqref{24}, \eqref{25}). These facts are used repeatedly in this paper. In Theorem 1, we show that, for $r=0$, the degenerate $r$-Whitney numbers reduce to the degenerate Stirling numbers. In Theorem 2, it is shown that, for $m=1$, the degenerate $r$-Whitney numbers reduce to the degenerate $r$-Stirling numbers. In Theorems 3 and 4, we express the degenerate $r$-Whitney numbers in terms of the degenerate Stirling numbers. In Theorems 6,7,9 and 10, we derive recurrence relations for the degenerate $r$-Whitney numbers. In Theorems 8 and 11, we prove some recurrence relations for the degenerate $r$-Stirling numbers. In Theorem 12, we deduce the orthogonality relations for the degenerate $r$-Whitney numbers from which the inversion relations are derived in Theorem 13. A certain expression is obtained for the degenerate $r$-Dowling polynomials in Theorem 14. An expression for the degenerate $r$-Whitney numbers of the second kind and a yet another recurrence relation for those numbers are given in Theorem 15 and Theorem 16, respectively. \par
For any $\lambda\in\mathbb{R}$, the degenerate exponential functions are defined by 
\begin{equation}
e_{\lambda}^{x}(t)=\sum_{k=0}^{\infty}\frac{(x)_{k,\lambda}}{k!}t^{k},\quad e_{\lambda}(t)=e_{\lambda}^{1}(t),\quad (\mathrm{see}\ [11]).\label{1}
\end{equation}
Here we recall that the generalized falling factorial sequence is defined by
\begin{equation}
(x)_{0,\lambda}=1,\quad (x)_{n,\lambda}=x(x-\lambda)(x-(n-1)\lambda),\quad (n\ge 1), \label{2}
\end{equation}
and that the generalized rising factorial sequence is given by 
\begin{equation}
\langle x\rangle_{0,\lambda}=1,\quad \langle x\rangle_{n,\lambda}=x(x+\lambda)(x+2\lambda)\cdots(x+(n+1)\lambda),\quad (n\ge 1). \label{3}	
\end{equation} \par
The degenerate Stirling numbers of the first kind are defined by 
\begin{equation}
(x)_{n}=\sum_{k=0}^{n}S_{1,\lambda}(n,k)(x)_{k,\lambda},\quad (n\ge 0),\quad (\mathrm{see}\ [1-16]),\label{4}
\end{equation}
where $(x)_{0}=1,\ (x)_{n}=x(x-1)\cdots(x-n+1),\ (n\ge 1)$. \\ 
Note that $\displaystyle\lim_{\lambda\rightarrow 0}S_{1,\lambda}(n,k)=S_{1}(n,k)\displaystyle$, where $S_{1}(n,k)$ are the ordinary Stirling number of the first kind given by 
\begin{displaymath}
(x)_{n}=\sum_{k=0}^{n}S_{1}(n,k)x^{k},\quad (n\ge 0),\quad (\mathrm{see}\ [8,10,12,13]). 
\end{displaymath}
As the inversion formula of \eqref{4}, the degenerate Stirling numbers of the second kind are defined by 
\begin{equation}
(x)_{n,\lambda}=\sum_{k=0}^{n}S_{2,\lambda}(n,k)(x)_{k},\quad (n\ge 0),\quad (\mathrm{see}\ [8]).\label{5}
\end{equation}
Note that $\displaystyle\lim_{\lambda\rightarrow 0}S_{2,\lambda}(n,k)=S_{2}(n,k)\displaystyle$, where $S_{2}(n,k)$ are the ordinary Stirling numbers of the second kind defined by 
\begin{displaymath}
x^{n}=\sum_{k=0}^{n}S_{2}(n,k)(x)_{k},\quad (\mathrm{see}\ [10,13]). 
\end{displaymath}
For $r\in\mathbb{N}\cup\{0\}$, the unsigned degenerate $r$-Stirling numbers of the first kind are defined by 
\begin{equation}
\langle x+r\rangle_{n}=\sum_{k=0}^{n}{n+r \brack k+r}_{r,\lambda}\langle x\rangle_{k,\lambda},\quad (n\ge 0),\quad (\mathrm{see}\ [13,14]),\label{6}
\end{equation}
where $\langle x\rangle_{0}=1,\ \langle x\rangle_{n}=x(x+1)\cdots(x+n-1),\ (n\ge 1)$. \\
In view of \eqref{5}, the degenerate $r$-Stirling numbers of the second kind are defined by 
\begin{equation}
(x+r)_{n,\lambda}=\sum_{k=0}^{n}{n+r \brace k+r}_{r,\lambda}(x)_{k},\quad (n\ge 0),\quad (\mathrm{see}\ [13,14]). \label{7}	
\end{equation} \par
The Whitney numbers $W_{m}(n,k)$ and $V_{m}(n,k)$ satisfy the following Stirling number-like relations: 
\begin{equation}
(mx+1)^{n}=\sum_{k=0}^{n}W_{m}(n,k)m^{k}(x)_{k},\quad (n\ge 0),\label{8}
\end{equation}
and 
\begin{equation}
m^{n}(x)_{n}=\sum_{k=0}^{n}V_{m}(n,k)(mx+1)^{k},\quad (\mathrm{see}\ [2,3,4,5,10]).\label{9}
\end{equation}
For $n\ge 0$, Dowling polynomials are defined by 
\begin{equation}
D_{m}(n,x)=\sum_{k=0}^{n}W_{m}(n,k)x^{k},\quad (n\ge 0),\quad (\mathrm{see}\ [2,3,4,5,10]).\label{10}
\end{equation}
We remark that, as generalizations of Whitney numbers of both kinds, the $r$-Whitney numbers of the second $W_{m}^{(r)}(n,k)$ and of the first kind $V_{m}^{(r)}(n,k)$ were investigated by several authors (see [2,3,4,5,10]). They are defined by
\begin{equation}
(mx+r)^{n}=\sum_{k=0}^{n}W_{m}^{(r)}(n,k)m^{k}(x)_{k},\quad (\mathrm{see}\ [10]),\label{11}
\end{equation}
and 
\begin{equation} 
m^{n}(x)_{n}=\sum_{k=0}^{n}V_{m}^{(r)}(n,k)(mx+r)^{k},\quad (n\ge 0),\label{12}
\end{equation}
where $r$ is a non-negative integer. \par 
Recently, the degenerate Whitney numbers of the second kind $W_{m,\lambda}(n,k)$ and of the first kind $V_{m,\lambda}(n,k)$ are introduced in [10] as degenerate versions of $W_{m}(n,k)$ and $V_{m}(n,k)$, respectively. Indeed, they are defined by
\begin{align}
(mx+1)_{n,\lambda}=\sum_{k=0}^{n}W_{m,\lambda}(n,k)m^{k}(x)_{k},\quad(n \ge 0), \label{13}
\end{align}
and
\begin{align}
m^{n}(x)_{n}=\sum_{k=0}^{n}V_{m,\lambda}(n,k)(mx+1)_{k,\lambda},\quad (\mathrm{see}\ [10])\label{14}.
\end{align}
As further generalizations of the degenerate Whitney numbers of both kinds, the degenerate $r$-Whitney numbers of the second kind $W_{m,\lambda}^{(r)}(n,k)$ and of the first kind $V_{m,\lambda}^{(r)}(n,k)$ are defined by
\begin{equation}
(mx+r)_{n,\lambda}=\sum_{k=0}^{n}W_{m,\lambda}^{(r)}(n,k)m^{k}(x)_{k},\quad (\mathrm{see}\ [10]),\label{15}
\end{equation}
and 
\begin{equation} 
m^{n}(x)_{n}=\sum_{k=0}^{n}V_{m,\lambda}^{(r)}(n,k)(mx+r)_{k,\lambda},\quad (n\ge 0),\label{16}
\end{equation}
where $r$ is a non-negative integer. \\ 
Note that $\displaystyle\lim_{\lambda\rightarrow 0}V_{m,\lambda}^{(r)}(n,k)=V_{m}^{(r)}(n,k)\displaystyle$ and $\displaystyle\lim_{\lambda\rightarrow 0}W_{m,\lambda}^{(r)}(n,k)=W_{m}^{(r)}(n,k)\displaystyle$. 
In [10], the degenerate Dowling polynomials are defined by 
\begin{equation}
D_{m,\lambda}(n,x)=\sum_{k=0}^{n}W_{m,\lambda}(n,k)x^{k},\quad (n\ge 0). \label{17}	
\end{equation}
When $x=1$, $D_{m,\lambda}(n,1)=D_{m,\lambda}(x)$ are called the degenerate Dowling numbers. \par
The boson operators $a^{+}$ and $a$ satisfy the following commutation relation: 
\begin{equation}
[a,a^{+}]=aa^{+}-a^{+}a=1.\label{18}
\end{equation}
The normal ordering of an integral power of the number operator $a^{+}a$ in terms of boson operators $a$ and $a^{+}$ can be written in the form 
\begin{equation}
(a^{+}a)^{k}=\sum_{l=0}^{k}S_{2}(k,l)(a^{+})^{l}a^{l},\quad (\mathrm{see}\ [7,9,12,16]).\label{19}
\end{equation}
The number states $|m\rangle,\ m=0,1,2,\dots$, are defined as 
\begin{equation}
a|m\rangle=\sqrt{m}|m-1\rangle,\quad a^{+}|m\rangle=\sqrt{m+1}|m+1\rangle.\label{20}
\end{equation}
The coherent states $|z\rangle$, where $z$ is a complex number, satisfy $a|z\rangle=z|z\rangle,\ \langle z|z\rangle=1$. To show a connection to coherent states, we recall that the harmonic oscillator has Hamiltonian $\hat{n}=a^{+}a$ (neglecting the zero point energy) and the usual eigenstates $|n\rangle$ (for $n\in\mathbb{N})$ satisfying $\hat{n}|n=n|n\rangle$ and $\langle m|n\rangle=\delta_{m,n}$ where $\delta_{m,n}$ is the Kronecker's symbol. \par 

\section{Degenerate $r$-Whitney numbers and $r$-Dowling polynomials via boson operators}
First, we recall that the standard bosonic commutation relations $aa^{+}-a^{+}a=1$ can be considered formally, in a suitable space of functions $f$, by letting $a=\frac{d}{dx}$ and $a^{+}=x$ (the operator of multiplication by $x$). From the definition of coherent state $a|z\rangle=z|z\rangle$, we note that $\langle z|a=\langle z|\bar{z}$, where $z\in\mathbb{C}$. \par 
Let $f$ be a polynomial. Then we easily get 
\begin{align}
&\bigg(x\frac{d}{dx}\bigg)_{n,\lambda}f(x)=\sum_{k=0}^{n}S_{2,\lambda}(n,k)x^{k}\bigg(\frac{d}{dx}\bigg)^{k}f(x)\label{21} \\
&x^{n}\bigg(\frac{d}{dx}\bigg)^{n}f(x)=\sum_{k=0}^{n}S_{1,\lambda}(n,k)\bigg(x\frac{d}{dx}\bigg)_{k,\lambda}f(x),\quad (n\ge 0). \nonumber
\end{align}
In view of \eqref{21}, the normal ordering of the degenerate $n$-th power of the number operator $a^{+}a$ in terms of boson operators $a,a^{+}$ can be written in the form 
\begin{align}
&(a^{+}a)_{n,\lambda}=\sum_{k=0}^{n}S_{2,\lambda}(n,k)(a^{+})^{k}a^{k}, \label{22} \\
& (a^{+})^{n}a^{n}=\sum_{k=0}^{n}S_{1,\lambda}(n,k)(a^{+}a)_{k,\lambda}, \quad (n\ge 0). \label{23}
\end{align}
Now, \eqref{15} and \eqref{16} can be written as 
\begin{equation}
m^{n}(a^{+})^{n}a^{n}=\sum_{k=0}^{n}V_{m,\lambda}^{(r)}(n,k)(ma^{+}a+r)_{k,\lambda}, \label{24}
\end{equation}
and 
\begin{equation}
(ma^{+}a+r)_{n,\lambda}=\sum_{k=0}^{n}W_{m,\lambda}^{(r)}(n,k)m^{k} (a^{+})^{k}a^{k},\quad (n\ge 0).\label{25}
\end{equation}
In [10], it is shown that 
\begin{displaymath}
	V_{m,\lambda}^{(r)}(n,k)=W_{m,\lambda}^{(r)}(n,k)=0\quad \textrm{if $n<k$ or $k<0$},  
\end{displaymath}
and $V_{m,\lambda}^{(r)}(0,0)=W_{m,\lambda}^{(r)}(0,0)=1$. \par 
For $k\in\mathbb{N}$, by \eqref{20}, we get 
\begin{equation}
(a^{+})^{k}a^{k}|m\rangle=(m)_{k}|m\rangle,\label{26}
\end{equation}
where $m$ is nonnegative integer. \par 
From \eqref{24} and \eqref{25}, we note that 
\begin{align}
&\sum_{k=0}^{n}V_{m,\lambda}^{(0)}(n,k)(ma^{+}a)_{k,\lambda}=m^{n}(a^{+})^{n}a^{n} \label{27} \\
&=m^{n}\sum_{k=0}^{n}S_{1,\frac{\lambda}{m}}(n,k)(a^{+}a)_{k,\frac{\lambda}{m}}=\sum_{k=0}^{n}m^{n-k}S_{1,\frac{\lambda}{m}}(n,k)(ma^{+}a)_{k,\lambda},\nonumber
\end{align}
and 
\begin{align}
\sum_{k=0}^{n}W_{m,\lambda}^{(0)}(n,k)m^{k}(a^{+})^{k}a^{k}&=(ma^{+}a)_{m,\lambda}=m^{n}(a^{+}a)_{n,\frac{\lambda}{m}}\label{28} \\
&=m^{n}\sum_{k=0}^{n}S_{2,\frac{\lambda}{m}}(n,k)(a^{+})^{k}a^{k}. \nonumber
\end{align}
Therefore, by \eqref{27} and \eqref{28}, we obtain the following theorem.
\begin{theorem}
For $n\ge 0$, we have 
\begin{displaymath}
V_{m,\lambda}^{(0)}(n,k)=m^{n-k}S_{1,\frac{\lambda}{m}}(n,k),\quad W_{m,\lambda}^{(0)}(n,k)=m^{n-k}S_{2,\frac{\lambda}{m}}(n,k).
\end{displaymath}	
\end{theorem}
 From \eqref{6}, we note that 
 \begin{align}
 (x-r)_{n}&=(-1)^{n}\langle -x+r\rangle_{n}=(-1)^{n}\sum_{k=0}^{n}{n+r \brack k+r}_{r,\lambda}\langle -x\rangle_{k,\lambda}\label{29}\\
 &=\sum_{k=0}^{n}(-1)^{n-k}{n+r \brack k+r}_{r,\lambda}(x)_{k,\lambda},	\nonumber
 \end{align}
 and hence that
\begin{align}
(x)_{n}=\sum_{k=0}^{n}(-1)^{n-k}{n+r \brack k+r}_{r,\lambda}(x+r)_{k,\lambda}. \label{30}
\end{align}
Let $f$ be a polynomial and let $D=\frac{d}{dx}$. Then, by \eqref{30}, we get 
\begin{align}
x^{n}D^{n}f(x)=\sum_{k=0}^{n}(-1)^{n-k}{n+r \brack k+r}_{r,\lambda}(xD+r)_{k,\lambda}f(x), \label{31}
\end{align}
and, by inversion, we have
\begin{align}
(xD+r)_{n,\lambda}f(x)=\sum_{k=0}^{n}{n+r \brace k+r}_{r,\lambda}x^{k}D^{k}f(x).\label{32}
\end{align}
From \eqref{31} and \eqref{32}, we note that 
\begin{align}
(a^{+})a^{n}=\sum_{k=0}^{n}(-1)^{n-k}{n+r \brack k+r}_{r,\lambda}(a^{+}a+r)_{k,\lambda}, \label{33}
\end{align}
and 
\begin{align}
(a^{+}a+r)_{n,\lambda}=\sum_{k=0}^{n}{n+r \brace k+r}_{r,\lambda}(a^{+})^{k}a^{k}. \label{34}
\end{align}
By \eqref{24}, \eqref{25}, \eqref{33} and \eqref{34}, we get 
\begin{align}
(a^{+})^{n}a^{n}&=\sum_{k=0}^{n}V_{1,\lambda}^{(r)}(n,k)(a^{+}a+r)_{k,\lambda} \label{35}\\
&=\sum_{k=0}^{n}(-1)^{n-k}{n+r \brack k+r}_{r,\lambda}(a^{+}a+r)_{k,\lambda}, \nonumber
\end{align}
and 
\begin{align}
(a^{+}a+r)_{n,\lambda}&=\sum_{k=0}^{n}W_{1,\lambda}^{(r)}(n,k)(a^{+})^{k}a^{k} \label{36} \\
&=\sum_{k=0}^{n}{n+r \brace k+r}_{r,\lambda}(a^{+})^{k}a^{k}.\nonumber
\end{align}
Therefore, by \eqref{35} and \eqref{36}, we obtain the following theorem. 
\begin{theorem}
For $n,k\ge 0$ with $n\ge k$, we have 
\begin{align*}
	&V_{1,\lambda}^{(r)}(n,k)=(-1)^{n-k}{n+r \brack k+r}_{r,\lambda}, \\
	&W_{1,\lambda}^{(r)}(n,k)={n+r \brace k+r}_{r,\lambda}. 
\end{align*}
\end{theorem}
From \eqref{22} and \eqref{24}, we have 
\begin{align}
&\sum_{k=0}^{n}V_{m,m\lambda}^{(r)}(n,k)(ma^{+}a+r)_{k,m\lambda}=m^{n}(a^{+})^{n}a^{n}\label{37}\\ 
 &=m^{n}\sum_{i=0}^{n}S_{1,\lambda}(n,i)(a^{+}a)_{i,\lambda}=m^{n}\sum_{i=0}^{n}S_{1,\lambda}(n,i)\bigg(\frac{ma^{+}a+r-r}{m}\bigg)_{i,\lambda}\nonumber\\
&=m^{n}\sum_{i=0}^{n}S_{1,\lambda}(n,i)\bigg(\frac{1}{m}\bigg)^{i}(ma^{+}a+r-r)_{i,m\lambda} \nonumber \\
&=\sum_{i=0}^{n}S_{1,\lambda}(n,i) m^{n-i}\sum_{k=0}^{i}\binom{i}{k}(ma^{+}a+r)_{k,m\lambda}(-r)_{i-k,m\lambda}\nonumber\\
&=\sum_{k=0}^{n}\bigg(\sum_{i=k}^{n}S_{1,\lambda}(n,i)m^{n-i}\binom{i}{k}(-r)_{i-k,m\lambda}\bigg)(ma^{+}a+r)_{k,m\lambda}.\nonumber 
\end{align}
Therefore, by comparing the coefficients on both sides of \eqref{37}, we obtain the following theorem. 
\begin{theorem}
\begin{displaymath}
V_{m,m\lambda}^{(r)}(n,k)=\sum_{i=k}^{n}\binom{i}{k}m^{n-i}S_{1,\lambda}(n,i)(-r)_{i-k,m\lambda}.
\end{displaymath}
\end{theorem}
By \eqref{23} and \eqref{25}, we get 
\begin{align}
&\sum_{k=0}^{n}W_{m,\lambda}^{(r)}(n,k)m^{k}(a^{+})^{k}a^{k}=(ma^{+}a+r)_{n,\lambda}\label{38}\\
&=\sum_{i=0}^{n}\binom{n}{i}(r)_{n-i,\lambda}(ma^{+}a)_{i,\lambda}=\sum_{i=0}^{n}\binom{n}{i}(r)_{n-i,\lambda}m^{i}(a^{+}a)_{i,\frac{\lambda}{m}}\nonumber \\
&=\sum_{i=0}^{n}\binom{n}{i}(r)_{n-i,\lambda}m^{i}\sum_{k=0}^{i}S_{2,\frac{\lambda}{m}}(i,k)(a^{+})^{k}a^{k} \nonumber \\
&=\sum_{k=0}^{n}\bigg\{\sum_{i=k}^{n}\binom{n}{i}(r)_{n-i,\lambda}m^{i}S_{2,\frac{\lambda}{m}}(i,k)\bigg\}(a^{+})^{k}a^{k}. \nonumber
\end{align}
Therefore, by comparing the coefficients on both sides of \eqref{38}, we obtain the following theorem. 
\begin{theorem}
For $n,k\ge 0$ with $n\ge k$, we have 
\begin{displaymath}
W_{m,\lambda}^{(r)}(n,k)=\sum_{i=k}^{n}\binom{n}{i}(r)_{n-i,\lambda}m^{i-k}S_{2,\frac{\lambda}{m}}(i,k). 
\end{displaymath}	
\end{theorem}
The degenerate unsigned Stirling numbers of the first kind are defined by 
\begin{equation}
\langle x\rangle_{n}=\sum_{k=0}^{n}{n \brack k}_{\lambda}\langle x\rangle_{k,\lambda},\quad (n\ge 0).\label{39}
\end{equation}
Thus, by \eqref{4} and \eqref{39}, we get 
\begin{equation}
(-1)^{n-i}S_{1,\lambda}(n,i)={n \brack i}_{\lambda},\quad (n \ge i\ge 0). \label{40}
\end{equation}
Let us take $m=1$. Then, by Theorems 2, 3 and 4, we get 
\begin{equation}
\sum_{i=k}^{n}\binom{i}{k}S_{1,\lambda}(n,i)(-r)_{i-k,\lambda}=V_{1,\lambda}^{(r)}(n,k)=(-1)^{n-k}{n+r \brack k+r}_{r,\lambda}, \label{41}
\end{equation}
and 
\begin{equation}
\sum_{i=k}^{n}\binom{n}{i}(r)_{n-i,\lambda}S_{2,\lambda}(i,k)=W_{1,\lambda}^{(r)}(n,k)={n+r \brace k+r}_{r,\lambda}.\label{42}
\end{equation}
Therefore, by \eqref{40}, \eqref{41} and \eqref{42}, we obtain the following theorem. 
\begin{theorem}
For $n,k\ge 0$ with $n\ge k$, we have 
\begin{displaymath}
{n+r \brack k+r}_{r,\lambda}=\sum_{i=k}^{n}(-1)^{i-k}(-r)_{i-k,\lambda}{n \brack i}_{\lambda}\binom{i}{k},
\end{displaymath}
and 
\begin{displaymath}
{n+r \brace k+r}_{r,\lambda}=\sum_{i=k}^{n}\binom{n}{i}(r)_{n-i,\lambda}S_{2,\lambda}(i,k).
\end{displaymath}
\end{theorem}
We claim that
\begin{equation}
[a^{n},a^{+}]=a^{n}a^{+}-a^{+}a^{n}=na^{n-1},\quad (n\in\mathbb{N}).\label{43}	
\end{equation}
Note that $n=1$ case is just the commutation relation. Assume that \eqref{43} holds true.\\
Then we have
\begin{align*}
[a^{n+1},a^{+}]&=a^{n+1}a^{+}-a^{+}a^{n+1}=a(a^{n}a^{+})-a^{+}a^{n+1}\\
&=a(a^{+}a^{n}+na^{n-1})-a^{+}a^{n+1}\\
&=na^{n}+(aa^{+}-a^{+}a)a^{n}=(n+1)a^{n}.
\end{align*}
In a similar manner to \eqref{43}, we can show that
\begin{equation}
[a,(a^{+})^{n}]=a(a^{+})^{n}-(a^{+})^{n}a=n(a^{+})^{n-1},\quad (n\in\mathbb{N}). \label{44}
\end{equation}
By \eqref{24} and \eqref{43}, we get 
\begin{align}
&\sum_{k=0}^{n+1}V_{m,\lambda}^{(r)}(n+1,k)(ma^{+}a+r)_{k,\lambda}=m^{n+1}(a^{+})^{n+1}a^{n+1} \label{45} \\
&=m^{n+1}(a^{+})^{n}(a^{+}a^{n})a=m^{n+1}(a^{+})^{n}(a^{n}a^{+}-na^{n-1})a	\nonumber \\
&= m^{n}(a^{+})^{n}a^{n}(ma^{+}a)-(mn)m^{n}(a^{+})^{n}a^{n} \nonumber \\
&=\sum_{k=0}^{n}V_{m,\lambda}^{(r)}(n,k)(ma^{+}a+r)_{k,\lambda}(ma^{+}a+r-k\lambda+k\lambda-r)\nonumber \\
&\quad -mn\sum_{k=0}^{n}V_{m,\lambda}^{(r)}(n,k)(ma^{+}a+r)_{k,\lambda} \nonumber \\
&=\sum_{k=0}^{n}V_{m,\lambda}^{(r)}(n,k)(ma^{+}a+r)_{k+1,\lambda}+\sum_{k=0}^{n}(k\lambda-r-nm)V_{m,\lambda}^{(r)}(n,k)(ma^{+}a+r)_{k,\lambda}\nonumber\\
&=\sum_{k=1}^{n+1}V_{m,\lambda}^{(r)}(n,k-1)(ma^{+}a+r)_{k,\lambda}+\sum_{k=0}^{n}(k\lambda-r-nm)V_{m,\lambda}^{(r)}(n,k)(ma^{+}a+r)_{k,\lambda}\nonumber\\
&=\sum_{k=0}^{n+1}\Big\{V_{m,\lambda}^{(r)}(n,k-1)+(k\lambda-r-mn)V_{m,\lambda}^{(r)}(n,k)\big\}(ma^{+}a+r)_{k,\lambda}.\nonumber
\end{align}
Therefore, by comparing the coefficients on both sides of \eqref{45}, we obtain the following theorem. 
\begin{theorem}
For $n,k\in\mathbb{N}$ with $n\ge k$, we have 
\begin{displaymath}
V_{m,\lambda}^{(r)}(n+1,k)=V_{m,\lambda}^{(r)}(n,k-1)+(k\lambda-r-mn)V_{m,\lambda}^{(r)}(n,k).
\end{displaymath}	
\end{theorem}
From \eqref{25}, we note that 
\begin{align}
&\sum_{k=0}^{n+1}m^{k}W_{m,\lambda}^{(r)}(n+1,k)(a^{+})^{k}a^{k}=(ma^{+}a+r)_{n+1,\lambda}=(ma^{+}a+r)_{n,\lambda}(ma^{+}a+r-n\lambda) \label{46} \\
&=\sum_{k=0}^{n}m^{k}W_{m,\lambda}^{(r)}(n,k)(a^{+})^{k}a^{k}ma^{+}a+(r-n\lambda)\sum_{k=0}^{n}m^{k}W_{m,\lambda}^{(r)}(n,k)(a^{+})^{k}a^{k} \nonumber \\
&=\sum_{k=0}^{n}m^{k+1}W_{m,\lambda}^{(r)}(n,k)(a^{+})^{k}(a^{k}a^{+})a+(r-n\lambda)\sum_{k=0}^{n}m^{k}W_{m,\lambda}^{(r)}(n,k)(a^{+})^{k}a^{k}\nonumber \\
&=\sum_{k=0}^{n}m^{k+1}W_{m,\lambda}^{(r)}(n,k)(a^{+})^{k}(a^{+}a^{k}+ka^{k-1})a+(r-n\lambda)\sum_{k=0}^{n}m^{k}W_{m,\lambda}^{(r)}(n,k)(a^{+})^{k}a^{k}\nonumber \\
&=\sum_{k=0}^{n}m^{k+1}W_{m,\lambda}^{(r)}(n,k)(a^{+})^{k+1}a^{k+1}+\sum_{k=0}^{n}m^{k}W_{m,\lambda}^{(r)}(n,k)(mk+r-n\lambda)(a^{+})^{k}a^{k}\nonumber \\
&=\sum_{k=1}^{n+1}m^{k}W_{m,\lambda}^{(r)}(n,k-1)(a^{+})^{k}a^{k}+\sum_{k=0}^{n}m^{k}W_{m,\lambda}^{(r)}(n,k)(mk+r-n\lambda)(a^{+})^{k}a^{k} \nonumber \\
&=\sum_{k=0}^{n+1}m^{k}\Big\{W_{m,\lambda}^{(r)}(n,k-1)+(mk+r-n\lambda)W_{n,\lambda}^{(r)}(n,k)\Big\}(a^{+})^{k}a^{k}. \nonumber
\end{align}
Therefore, by comparing the coefficients on both sides of \eqref{46}, we obtain the following theorem. 
\begin{theorem}
For $n,k\in\mathbb{N}$ with $n\ge k$, we have 
\begin{displaymath}
W_{m,\lambda}^{(r)}(n+1,k)=W_{m,\lambda}^{(r)}(n,k-1)+(mk+r-n\lambda)W_{m,\lambda}^{(r)}(n,k).
\end{displaymath}	
\end{theorem}
From Theorem 6, we note that 
\begin{align*}
V_{m,\lambda}^{(r)}(n+1,0)&=-(r+mn)V_{m,\lambda}^{(r)}(n,0)=(-1)^{2}(r+mn)(r+m(n-1))V_{m,\lambda}^{(r)}(n-1,0)\\
&=\cdots=(-1)^{n+1}(r+mn)(r+m(n-1))\cdots rV_{m,\lambda}^{(r)}(0,0).	
\end{align*}
Thus, we get 
\begin{equation}
V_{m,\lambda}^{(r)}(n,0)=(-1)^{n}\prod_{k=0}^{n-1}(mk+r),\quad (n\in\mathbb{N}).\label{47}
\end{equation}
From Theorem 3 and \eqref{47}, we have 
\begin{align}
V_{m,\lambda}^{(r)}(n,0)&=\sum_{i=0}^{n}S_{1,\lambda/m}(n,i)m^{n-i}(-r)_{i,\lambda} \nonumber\\
&=(-1)^{n}\prod_{k=0}^{n-1}(mk+r),\quad (n\in\mathbb{N}).\label{48}	
\end{align}
By Theorem 7, we get 
\begin{align}
W_{m,\lambda}^{(r)}(n+1,0)&=(r-n\lambda)W_{m,\lambda}^{(r)}(n,0)=(r-n\lambda)(r-(n-1)\lambda)W_{n,\lambda}^{(r)}(n-1,0)  \label{49} \\
&=\cdots=(r-n\lambda)(r-(n-1)\lambda)\cdots rW_{m,\lambda}^{(r)}(0,0) \nonumber \\
&=(r)_{n+1,\lambda}. \nonumber	
\end{align}
Thus, we have 
\begin{equation}
W_{m,\lambda}^{(r)}(n,0)=(r)_{n,\lambda},\quad (n\in\mathbb{N}\cup\{0\}). \label{50}
\end{equation}
From Theorems 2 and 6, we have 
\begin{align}
(-1)^{n+1-k}{n+1+r \brack k+r}_{r,\lambda}&=V_{1,\lambda}^{(r)}(n+1,k)=V_{1,\lambda}^{(r)}(n,k-1)+(k\lambda-r-n)V_{1,\lambda}^{(r)}(n,k)\label{51}\\
&=(-1)^{n-k+1}{n+r \brack k+r-1}_{r,\lambda}+(k\lambda-r-n)(-1)^{n-k}{n+r \brack k+r}_{r,\lambda}.\nonumber 
\end{align}
Thus, we have 
\begin{equation}
{n+1+r\brack k+r}_{r,\lambda}={n+r \brack k+1+r}_{r,\lambda}+(r+n-k\lambda){n+r \brack k+r}_{r,\lambda}.\label{52}	
\end{equation}
Also, from Theorems 2 and 7, we get
\begin{align}
{n+r+1 \brace k+r}_{r,\lambda}&=W_{1,\lambda}^{(r)}(n+1,k)=W_{1,\lambda}^{(r)}(n,k-1)+(k+r-n\lambda)W_{1,\lambda}^{(r)}(n,k) \label{53}\\
&={n+r \brace k-1+r}_{r,\lambda}+(k+r-n\lambda){n+r \brace k+r}_{r,\lambda}.\nonumber	
\end{align}
Therefore, by \eqref{52} and \eqref{53}, we obtain the following theorem. 
\begin{theorem}
For $n,k\ge 0$ with $n\ge k$, we have 
\begin{align*}
	{n+r+1 \brack k+r}_{r,\lambda}={n+r \brack k+r-1}_{r,\lambda}+(r+n-k\lambda){n+r \brack k+r}_{r,\lambda},
\end{align*}	
and 
\begin{align*}
	{n+r+1 \brace k+r}_{r,\lambda}={n+r \brace k+r-1}_{r,\lambda}+(r+k-n\lambda){n+r \brace k+r}_{r,\lambda}.
\end{align*}
\end{theorem}
On the one hand, from \eqref{24} we have
\begin{align}
m^{n}(a^{+})^{n}a^{n} &=\sum_{k=0}^{n}V_{m,\lambda}^{(r)}(n,k)(ma^{+}a+r)_{k,\lambda} \label{54} \\
&=\sum_{k=0}^{n}V_{m,\lambda}^{(r)}(n,k)(ma^{+}a+r+1-1)_{k,\lambda} 	\nonumber \\
&=\sum_{k=0}^{n}V_{m,\lambda}^{(r)}(n,k)\sum_{l=0}^{k}\binom{k}{l}(-1)^{k-l}\langle 1\rangle_{k-l,\lambda}(ma^{+}a+r+1)_{l,\lambda} \nonumber \\
&=\sum_{l=0}^{n}(ma^{+}a+r+1)_{l,\lambda}\bigg(\sum_{k=l}^{n}\binom{k}{l}(-1)^{k-l}\langle 1\rangle_{k-l,\lambda}V_{m,\lambda}^{(r)})(n,k)\bigg)\nonumber.
\end{align}
On the other hand, from \eqref{24}, we get 
\begin{equation}
m^{n}(a^{+})^{n}a^{n}=\sum_{l=0}^{n}V_{m,\lambda}^{(r+1)}(n,l)(ma^{+}a+r+1)_{l,\lambda}. \label{55}
\end{equation}
Therefore, by \eqref{54} and \eqref{55}, we obtain the following theorem. 
\begin{theorem}
For $n,l\in\mathbb{Z}$ with $n\ge l\ge 0$, we have 
\begin{displaymath}
V_{m,\lambda}^{(r+1)}(n,l)=\sum_{k=l}^{n}\binom{k}{l}(-1)^{k-l}\langle 1\rangle_{k-l,\lambda}V_{m,\lambda}^{(r)}(n,k). 
\end{displaymath}	
\end{theorem}
From \eqref{25}, we note that 
\begin{align}
&\sum_{k=0}^{n}m^{k}W_{m,\lambda}^{(r+1)}(n,k)(a^{+})^{k}a^{k}=(ma^{+}a+r+1)_{n,\lambda}\label{56} \\
&=\sum_{l=0}^{n}\binom{n}{l}(1)_{n-l,\lambda}(ma^{+}a+r)_{l,\lambda}=\sum_{l=0}^{n}\binom{n}{l}(1)_{n-l,\lambda}\sum_{k=0}^{l}m^{k}W_{m,\lambda}^{(r)}(l,k)(a^{+})^{k}a^{k} \nonumber\\
&=\sum_{k=0}^{n}m^{k}(a^{+})^{k}a^{k}\bigg(\sum_{l=k}^{n}\binom{n}{l}(1)_{n-l,\lambda}W_{m,\lambda}^{(r)}(l,k)\bigg).\nonumber
\end{align}
Therefore, by comparing the coefficients on both sides of \eqref{56}, we obtain the following theorem.
\begin{theorem}
For $n,k\in\mathbb{Z}$ with $n\ge k\ge 0$, we have 
\begin{displaymath}
	W_{m,\lambda}^{(r+1)}(n,k)=\sum_{l=k}^{n}\binom{n}{l}(1)_{n-l,\lambda}W_{m,\lambda}^{(r)}(l,k).
\end{displaymath}	
\end{theorem}
From Theorems 2 and 9, we note that 
\begin{align}
(-1)^{n-l}{n+r+1 \brack l+r+1}_{r+1,\lambda}&=V_{1,\lambda}^{(r+1)}(n,l)=\sum_{k=l}^{n}\binom{k}{l}(-1)^{k-l}\langle 1\rangle_{k-l,\lambda}V_{1,\lambda}^{(r)}(n,k) \nonumber \\
&=\sum_{k=l}^{n}\binom{k}{l}(-1)^{k-l}\langle 1\rangle_{k-l,\lambda}(-1)^{n-k}{n+r \brack k+r}_{r,\lambda} \label{57}\\
&=\sum_{k=l}^{n}\binom{k}{l}(-1)^{n-l}\langle 1\rangle_{k-l,\lambda}{n+r \brack k+r}_{r,\lambda},\nonumber
\end{align}
and from Theorems 2 and 10, we observe that
\begin{align}
{n+r+1 \brace k+r+1}_{r+1,\lambda}=W_{1,\lambda}^{(r+1)}(n,k)&=\sum_{l=k}^{n}\binom{n}{l}W_{1,\lambda}^{(r)}(l,k)(1)_{n-l,\lambda} \label{58} \\
&=\sum_{l=k}^{n}\binom{n}{l}{l+r \brace k+r}_{r,\lambda}(1)_{n-l,\lambda}.\nonumber
\end{align}
Therefore, by \eqref{57} and \eqref{58}, we obtain the following theorem. 
\begin{theorem}
For $n,l\in\mathbb{Z}$ with $n\ge l\ge 0$, we have 
\begin{displaymath}
{n+r+1 \brack l+r+1}_{r+1,\lambda}=\sum_{k=l}^{n}\binom{k}{l}\langle 1\rangle_{k-l,\lambda}{n+r \brack k+r}_{r,\lambda},
\end{displaymath}
and 
\begin{displaymath}
	{n+r+1 \brace l+r+1}_{r+1,\lambda}=\sum_{j=l}^{n}\binom{n}{j}{j+r \brace l+r}_{r,\lambda}(1)_{n-j,\lambda}.
\end{displaymath}
\end{theorem}
On the one hand, from \eqref{24} and \eqref{25}, we have 
\begin{align}
(ma^{+}a+r)_{n,\lambda}&=\sum_{k=0}^{n}W_{m,\lambda}^{(r)}(n,k)m^{k}(a^{+})^{k}a^{k} \label{59} \\
&=\sum_{k=0}^{n}W_{m,\lambda}^{(r)}(n,k)\sum_{j=0}^{k}V_{m,\lambda}^{(r)}(k,j)(ma^{+}a+r)_{j,\lambda}\nonumber \\
&=\sum_{j=0}^{n}\bigg(\sum_{k=j}^{n}W_{m,\lambda}^{(r)}(n,k)V_{m,\lambda}^{(k)}(k,j)\bigg)(ma^{+}a+r)_{j,\lambda}.\nonumber	
\end{align}
On the other hand, by the same method as \eqref{59}, we get 
\begin{align}
m^{n}(a^{+})^{n}a^{n} &= \sum_{k=0}^{n}V_{m,\lambda}^{(r)}(n,k)(ma^{+}a+r)_{k,\lambda} \label{60} \\
&=\sum_{k=0}^{n}V_{m,\lambda}^{(r)}(n,k)\sum_{j=0}^{k}	W_{m,\lambda}^{(r)}(k,j)m^{j}(a^{+})^{j}a^{j} \nonumber \\
&=\sum_{j=0}^{n}m^{j}\bigg(\sum_{k=j}^{n}V_{m,\lambda}^{(r)}(n,k)W_{m,\lambda}^{(r)}(k,j)\bigg)(a^{+})^{j}a^{j}.\nonumber
\end{align}
Therefore, by \eqref{59} and \eqref{60}, we obtain the following theorem. 
\begin{theorem}
For $n,j\ge 0$ with $n \ge j$, we have 
\begin{displaymath}
\sum_{k=j}^{n}W_{m,\lambda}^{(r)}(n,k)V_{m,\lambda}^{(j)}(k,j)=\delta_{n,j},\quad \sum_{k=j}^{n}V_{m,\lambda}^{(r)}(n,k)W_{m,\lambda}^{(r)}(k,j)=\delta_{n,j},
\end{displaymath}	
where $\delta_{n,j}$ is the Kronecker's symbol. 
\end{theorem}
Let $\displaystyle f_{n}=\sum_{l=0}^{n}V_{m,\lambda}^{(r)}(n,l)g_{l}\displaystyle$. Then, by Theorem 12, we get 
\begin{align}
\sum_{l=0}^{n}W_{m,\lambda}^{(r)}(n,l)f_{l}&=\sum_{l=0}^{n}W_{m,\lambda}^{(r)}(n,l)\sum_{j=0}^{l}V_{m,\lambda}^{(r)}(l,j)g_{j} \nonumber \\
&=\sum_{j=0}^{n}g_{j}\bigg(\sum_{l=j}^{n}W_{m,\lambda}^{(r)}(n,l)V_{m,\lambda}^{(r)}(l,j)\bigg)=g_{n}.\label{61}
\end{align}
Assume that $\displaystyle g_{n}=\sum_{l=0}^{n}W_{m,\lambda}^{(r)}(n,l)f_{l}\displaystyle$. Then we have 
\begin{align}
\sum_{l=0}^{n}V_{m,\lambda}^{(r)}(n,l)g_{l}&=\sum_{l=0}^{n}V_{m,\lambda}^{(r)}(n,l)\sum_{j=0}^{l}W_{m,\lambda}^{(r)}(l,j)f_{j} \label{62} \\
&=\sum_{j=0}^{n}f_{j}\bigg(\sum_{l=j}^{n}V_{m,\lambda}^{(r)}(n,l)W_{m,\lambda}^{(r)}(l,j)\bigg)\nonumber \\
&=f_{n}. \nonumber	
\end{align}
Therefore, by \eqref{61} and \eqref{62}. we obtain the following theorem. 
\begin{theorem}
For $n\ge 0$, we have 
\begin{displaymath}
f_{n}=\sum_{l=0}^{n}V_{m,\lambda}^{(r)}(n,l)g_{l}\ \Longleftrightarrow\ g_{n}=\sum_{l=0}^{n}W_{m,\lambda}^{(r)}(n,l)f_{l}.
\end{displaymath}
\end{theorem}
In view of \eqref{17}, we consider the degenerate $r$-Dowling polynomials given by 
\begin{equation}
D_{m,\lambda}^{(r)}(n,x)=\sum_{k=0}^{n}W_{m,\lambda}^{(r)}(n,k)x^{k},\quad (n\ge 0). \label{63}
\end{equation}
When $x=1$, $D_{m,\lambda}^{(r)}(n,1)=D_{m,\lambda}^{(r)}(n)$ are called the degenerate $r$-Dowling numbers. \par 
We recall that the coherent state 
\begin{displaymath}
|\gamma\rangle=e^{-\frac{|d|^{2}}{2}}\sum_{k=0}^{\infty}\frac{\gamma^{k}}{\sqrt{k!}}|k\rangle,
\end{displaymath}
where $\gamma$ is an arbitrary complex constant $a|\gamma\rangle=\gamma|\gamma\rangle=1$ and $\langle \gamma|\gamma\rangle=1$. \par 
From \eqref{25}, we note that 
\begin{align}
\langle \gamma|(ma^{+}a+\gamma)_{n,\lambda}|\gamma\rangle &=\sum_{k=0}^{n}m^{k}W_{m,\lambda}^{(r)}(n,k)\langle \gamma|(a^{+})^{k}a^{k}|\gamma\rangle \nonumber \\
&=\sum_{k=0}^{n}m^{k}W_{m,\lambda}^{(r)}(n,k)\langle \gamma|(\overline{\gamma})^{k}\gamma^{k}|\gamma\rangle \label{64} \\
&=\sum_{k=0}^{n}m^{k}W_{m,\lambda}^{(r)}(n,k)|\gamma|^{2k}\langle \gamma|\gamma\rangle \nonumber \\
&=\sum_{k=0}^{n}m^{k}W_{m,\lambda}^{(r)}(n,k)|\gamma|^{2k}.\nonumber
\end{align}
The left hand side of \eqref{64} can be evaluated using coherent state by 
\begin{align}
\langle \gamma|(ma^{+}a+\gamma)_{m,\lambda}|\gamma\rangle &=\sum_{k,l=0}^{\infty}\frac{d^{k}(\overline{d})^{l}}{\sqrt{k!}\sqrt{l!}}\langle l|k\rangle (mk+r)_{n,\lambda}e^{-\frac{|d|^{2}}{2}}e^{-\frac{|d|^{2}}{2}}\label{65} \\
&=e^{-|d|^{2}}\sum_{k=0}^{\infty}\frac{|d|^{2k}}{k!}(mk+r)_{n,\lambda}.\nonumber
\end{align}
Combining \eqref{63}, \eqref{64} and \eqref{65}, we have
\begin{align}
D_{m,\lambda}^{(r)}(n,m|\gamma|^{2})=e^{-|\gamma|^{2}}\sum_{k=0}^{\infty}\frac{|\gamma|^{2k}}{k!}(mk+r)_{n,\lambda}.\label{66}
\end{align}
As $\gamma$ is an arbitrary complex number, we obtain the following theorem from \eqref{66}.
\begin{theorem}
For $n\ge 0$, we have 
\begin{displaymath}
D_{m,\lambda}^{(r)}(n,x)=e^{-\frac{x}{m}}\sum_{k=0}^{\infty}\bigg(\frac{x}{m}\bigg)^{k}\frac{(mk+r)_{n,\lambda}}{k!}.	
\end{displaymath}
In particular, for $x=1$, 
\begin{displaymath}
D_{m,\lambda}^{(r)}(n)=e^{-\frac{1}{m}}\sum_{k=0}^{\infty}\frac{(mk+r)_{n,\lambda}}{m^{k}k!}.
\end{displaymath}
\end{theorem}
From Theorem 14, we can derive the generating function of $r$-Dowling polynomials as in the following:
\begin{displaymath}
e_{\lambda}^{r}e^{\frac{x}{m}(e_{\lambda}^{m}(t)-1)}=\sum_{n=0}^{\infty}D_{m,\lambda}^{(r)}(n,x)\frac{t^{n}}{n!}. 	
\end{displaymath}
Taking $x=my$ in Theorem 14, we have
\begin{align}
&\sum_{k=0}^{n}W_{m,\lambda}^{(r)}(n,k)m^{k}y^{k}=D_{m,\lambda}^{(r)}(n,my)=e^{-y}\sum_{k=0}^{n}y^{k}\frac{(mk+r)_{n,\lambda}}{k!}\label{67} \\
&=\sum_{l=0}^{\infty}\frac{(-y)^{l}}{l!}\sum_{k=0}^{\infty}\frac{y^{k}(mk+r)_{n,\lambda}}{k!}=\sum_{i=0}^{\infty}\frac{y^{i}}{i!}\Big(\sum_{k=0}^{i}(-1)^{i-k}\binom{i}{k}(mk+r)_{n,\lambda}\Big).\nonumber
\end{align}
By comparing the coefficients on both sides of \eqref{67}, we get 
\begin{equation}
m^{i}W_{m,\lambda}^{(r)}(n,i)=\frac{1}{i!}\sum_{k=0}^{i}\binom{i}{k}(-1)^{i-k}(mk+r)_{n,\lambda}.\label{68}
\end{equation}
Therefore, by \eqref{68}, we obtain the following theorem. 
\begin{theorem}
For $n,i\ge 0$ with $n \ge i$, we have 
\begin{displaymath}
W_{m,\lambda}^{(r)}(n,i)=\frac{1}{m^{i}i!}\sum_{k=0}^{i}\binom{i}{k}(-1)^{i-k}(mk+r)_{n,\lambda}.
\end{displaymath}
\end{theorem}
From Theorem 15, we have 
\begin{align}
\sum_{n=0}^{\infty}W_{m,\lambda}^{(r)}(n,k)\frac{t^{n}}{n!}&=\frac{1}{m^{k}k!}\sum_{j=0}^{k}\binom{k}{j}(-1)^{k-j}\sum_{n=0}^{\infty}(mj+r)_{n,\lambda}\frac{t^{n}}{n!} \label{69}	\\
&=\frac{1}{m^{k}k!}\sum_{j=0}^{k}\binom{k}{j}(-1)^{k-j}e_{\lambda}^{mj+r}(t) \nonumber \\
&=\frac{1}{m^{k}k!}e_{\lambda}^{r}(t)\big(e_{\lambda}^{m}(t)-1\big)^{k}.\nonumber
\end{align}
Let $\triangle$ be the difference operator defined by $\triangle f(x)=f(x+1)-f(x)$. Then we have 
\begin{equation}
\triangle^{n}f(x)=\sum_{k=0}^{n}\binom{n}{k}(-1)^{n-k}f(x+k),\quad (n\ge 0). \label{70}
\end{equation}
By Theorem 15 and \eqref{70}, we get 
\begin{align}
W_{m,\lambda}^{(r)}(n,k)&=\frac{1}{m^{k}k!}\sum_{l=0}^{k}(-1)^{k-l}\binom{k}{l}(ml+r)_{n,\lambda}\label{71}\\
&=\frac{m^{n-k}}{k!}\sum_{l=0}^{k}(-1)^{k-l}\binom{k}{l}\bigg(l+\frac{r}{m}\bigg)_{n,\frac{\lambda}{m}}=\frac{m^{n-k}}{k} \triangle^{k}\bigg(\frac{r}{m}\bigg)_{n,\frac{\lambda}{m}}.\nonumber
\end{align}
Thus, we have 
\begin{displaymath}
W_{m,\lambda}^{(r)}(n,k)=\frac{m^{n-k}}{k!}\triangle^{k}\bigg(\frac{r}{m}\bigg)_{n,\frac{\lambda}{m}},\quad (n \ge k\ge 0).
\end{displaymath}
For any $k\in\mathbb{N}$, we observe that 
\begin{equation}
a^{+}a(a^{+}a)_{k,\lambda}=a^{+}(aa^{+})_{k,\lambda}a.\label{72}	
\end{equation}
On the one hand, by \eqref{72} and \eqref{25}, we get 
\begin{align}
&a^{+}\big(m(1+a^{+}a)+r	\big)_{n,\lambda}a=a^{+}(maa^{+}+r)_{n,\lambda}a\label{73}\\
&=a^{+}a(ma^{+}a+r)_{n,\lambda}=\frac{1}{m}(ma^{+}a+r-n\lambda+n\lambda-r)(ma^{+}a+r)_{n,\lambda}\nonumber \\
&=\frac{1}{m}(ma^{+}a+r)_{n+1,\lambda}-\frac{r-n\lambda}{m}(ma^{+}a+r)_{n,\lambda} \nonumber \\
&=\sum_{k=0}^{n+1}m^{k-1}(a^{+})^{k}a^{k}\Big(W_{m,\lambda}^{(r)}(n+1,k)-(r-n\lambda)W_{m,\lambda}^{(r)}(n,k)\Big).\nonumber
\end{align}
On the other hand, by \eqref{25} and using $(x+y)_{n,\lambda}=\sum_{k=0}^{n}\binom{n}{k}(x)_{k,\lambda}(y)_{n-k,\lambda}$, we get 
\begin{align}
&a^{+}\big(m(1+a^{+}a)+r\big)_{n,\lambda}a=a^{+}(ma^{+}a+r+m)_{n,\lambda}a\label{74}\\
&=a^{+}\sum_{l=0}^{n}\binom{n}{l}(m)_{n-l,\lambda}(ma^{+}a+r)_{l,\lambda}a\nonumber\\
&=\sum_{l=0}^{n}\binom{n}{l}(m)_{n-l,\lambda}\sum_{k=0}^{l}m^{k}W_{m,\lambda}^{(r)}(l,k)(a^{+})^{k+1}a^{k+1} \nonumber\\
&=\sum_{l=0}^{n}\binom{n}{l}(m)_{n-l,\lambda}\sum_{k=1}^{l+1}m^{k-1}W_{m,\lambda}^{(r)}(l,k-1)(a^{+})^{k}a^{k} \nonumber\\
&=\sum_{k=1}^{n+1}\bigg(\sum_{l=k-1}^{n}\binom{n}{l}(m)_{n-l,\lambda}W_{m,\lambda}^{(r)}(l,k-1)\bigg)m^{k-1}(a^{+})^{k}a^{k}.\nonumber
\end{align}
From \eqref{73} and \eqref{74}, we obtain the following theorem. 
\begin{theorem}
For $n,k\in\mathbb{N}$ with $n \ge k$, we have 
\begin{displaymath}
W_{m,\lambda}^{(r)}(n+1,k)-(r-n\lambda)W_{m,\lambda}^{(r)}(n,k)=\sum_{l=k-1}^{n}\binom{n}{l}(m)_{n-l,\lambda}W_{m,\lambda}^{(r)}(l,k-1).
\end{displaymath}
\end{theorem}

\section{Conclusion}
Explorations for degenerate versions of some special numbers and polynomials, which began by Carlitz's work on the degenerate Bernoulli and degenerate Euler polynomials, have recently regained interests of many mathematicians.
In this paper, we derived some properties, recurrence relations, orthogonality relations and several identities on the degenerate $r$-Whitney numbers of both kinds from the normal ordering of $(ma^{+}a+r)_{n,\lambda}$ in \eqref{24} and its inversion in \eqref{25}. The normal ordering is expressed with the help of the degenerate $r$-Whitney numbers of the second kind, while its inversion is represented by means of the degenerate $r$-Whitney numbers of the first kind. In addition, the degenerate $r$-Dowling polynomials are considered as a natural extension of the degenerate $r$-Whitney numbers of the second kind and their properties are investigated.
\par
It is one of our future projects to continue to explore various degenerate versions of many special polynomials and numbers by using such tools as mathematical physics, combinatorial methods, generating functions, umbral calculus techniques, $p$-adic analysis, differential equations, special functions, probability theory and analytic number theory.


\begin{thebibliography}{9}
\bibitem{1}
Ali, S. T.; Antoine, J.-P.; Gazeau, J.-P. \emph{Coherent states, wavelets, and their generalizations.} Second edition. Theoretical and Mathematical Physics. Springer, New York, 2014. xviii+577 pp.
\bibitem{2}
Benoumhani, M. \emph{On Whitney numbers of Dowling lattices.} Discrete Math. \textbf{159} (1996), no. 1-3, 13-33,
\bibitem{3}
Dowling, T. A. \emph{A class of geometric lattices based on finite groups.} J. Combinatorial Theory Ser. B 14 (1973), 61-86. 
\bibitem{4}
Dowling, T. A. Erratum: ``\emph{A class of geometric lattices based on finite groups}"  J. Combinatorial Theory Ser. B 15 (1973), 211. 
\bibitem{5}
Dowling, T. A.; Wilson, R. M. \emph{Whitney number inequalities for geometric lattices.} Proc. Amer. Math. Soc. \textbf{47} (1975), 504-512.
\bibitem{6}
Graham, R. L.; Knuth, D. E.; Patashnik, O. \emph{Concrete mathematics. A foundation for computer science.} Second edition. Addison-Wesley Publishing Company, Reading, MA, 1994. xiv+657 pp. ISBN: 0-201-55802-5
\bibitem{7}
Hecht, K. T. \emph{The vector coherent state method and its application to problems of higher symmetries.} Lecture Notes in Physics, 290. Springer-Verlag, Berlin, 1987. vi+154 pp.
\bibitem{8}
Kim, D. S.; Kim, T. \emph{A note on a new type of degenerate Bernoulli numbers,} Russ. J. Math. Phys. \textbf{27} (2020), no.2, 227-235.
\bibitem{9}
Kim, T.; Kim, D. S. \emph{Some identities involving degenerate Stirling numbers arising from normal ordering.} arXiv:2204.01252
\bibitem{10}
Kim, T.; Kim, D. S. \emph{A study on degenerate Whitney numbers of the first and second kinds of Dowling lattices.} arXiv:2103.08904 .
\bibitem{11}
 Kim, T.; Kim, D. S. \emph{On some degenerate differential and degenerate difference operators.} Russ. J. Math. Phys. \textbf{29} (2022), no. 1, 37-46.
\bibitem{12}
Kim, T.; Kim, D. S.; Kim, H. K. \emph{Normal ordering of degenerate integral powers of number operator and its applications.} arXiv:2204.02595
\bibitem{13}
Kim, T.; Kim, D. S.; Lee, H.; Park, J.-W. \emph{A note on degenerate $r$-Stirling numbers.} J. Inequal. Appl. 2020, Paper No. 225, 12 pp.
\bibitem{14}
Kim, T.; Yao, Y.; Kim, D. S.; Jang, G.-W. \emph{Degenerate $r$-Stirling numbers and $r$-Bell polynomials.} Russ. J. Math. Phys. \textbf{25} (2018), no. 1, 44-58.
\bibitem{15}
Maslov, V. P. \emph{Deterministic and indeterministic theory of coherent structures in turbulence.} V. A. Steklov Mathematical Institute Preprint, I. Acad. Sci. USSR Dep. Math., Moscow, 1986. 47 pp.
\bibitem{16}
Perelomov, A. \emph{Generalized coherent states and their applications.} Texts and Monographs in Physics. SpringerVerlag, Berlin, 1986. xii+320 pp.{}
\end{thebibliography}
\end{document}